\documentclass[12pt]{amsart}
\usepackage{amsmath,amstext,amsfonts,amsthm,amssymb}
\usepackage{eucal}
\usepackage{diagrams}

\diagramstyle[scriptlabels]

\renewcommand{\subsubsection}[1]{\addtocounter{subsubsection}{1}
{\ \\[3pt]\bf \thesubsubsection. \  #1} }

\swapnumbers

\newtheorem{thm}[subsection]{Theorem}

\newtheorem{prp}[subsection]{Proposition}

{  \theoremstyle{definition}
           \newtheorem{rem}[subsection]{Remark}
           
}
\newcommand{\cE}{\mathcal{E}}

\newcommand{\cO}{\mathcal{O}}

\newcommand{\cZ}{\mathcal{Z}}

\newcommand{\id}{\mathrm{id}}

\newcommand{\Char}{\mathrm{char}\ }
\newcommand{\Cone}{\mathrm{cone}\ }

\newcommand{\Init}{\mathrm{Init}\ }
\renewcommand{\arg}{\mathrm{arg}\ }
\newcommand{\Arg}{\mathrm{Arg}\ }
\newcommand{\Aut}{\mathrm{Aut}\ }

\newcommand{\Op}{\mathtt{Op}}
\newcommand{\Com}{\mathtt{COM}}

\begin{document}

\title[Models for operads]{Erratum to ``Homological algebra of homotopy 
algebras''}
\author{V. Hinich}
\address{Dept. of mathematics, University of Haifa, Mount Carmel, Haifa 31905
Israel}

\email{hinich@math.haifa.ac.il}
\begin{abstract}Theorem 6.1.1 of~\cite{haha} on the existence of a
model structure on the category of operads is not valid in
the generality  claimed. We present a counter-example (due to B.~Fresse)
and a corrected version of the theorem. 

\end{abstract}

\maketitle
\section{}
\subsection{}
In ~\cite{haha} we claimed the following result which turned out to be 
partly wrong (see the counterexample below).

 Let $k$ be a commutative ring. Let  $C(k)$ denote the category of complexes
over $k$ and $\Op(C(k))$ the category of operads in $C(k)$.

{\bf (Erroneous) Theorem 6.1.1} of~\cite{haha}.
{\em The category of operads \\
$\Op(C(k))$ admits a closed model category
structure with componentwise quasi-isomorphisms as weak equivalence
and componentwise epimorphisms as fibrations.}

\subsection{}
\label{half}
Recall that a morphism $\cO\rTo \cO'$ of operads is called a cofibration
if it can be presented as a retract of a direct limit
$$ \cO=\cO_0\rTo\cO_1\rTo\ldots$$
where $\cO_{n+1}$ is obtained from $\cO_n$ by joining a collection of free 
variables $x_i$ with prescribed values of the differentials $dx_i\in\cO_n$.

The axioms of model category include two types of factorization
for morphisms and two types of lifting properties. One half of these 
properties, namely, the lifting property for cofibrations with respect to 
the acyclic fibrations, as well as decomposition of any morphism into
a cofibration followed by an acyclic fibration, is very easy and is valid
with no restriction. In particular, any operad in $C(k)$ admits a cofibrant
resolution. The other half of the properties can be deduced if one
knows that the map ~(\ref{joinM}) below is a quasi-isomorphism.

\subsection{}

Let $\cO$ be an operad and let $M=\Cone(\id_k)[i]$ be a shifted cone
of $\id_k$. Let $F(M,n)$ denote the free operad generated by the collection
$M$ of $n$-ary operations. Consider the natural map 
\begin{equation}
\label{joinM}
i:\cO\to\cO\coprod F(M,n).
\end{equation}

In order to establish the existence of the model category structure mentioned
above, one has to check that ~(\ref{joinM}) is a quasi-isomorphism.

\section{}
The following counterexample to ~\cite{haha}, Theorem 6.1.1, due to
 B.~Fresse, has been communicated to the author by C.~Berger.

Let $\Com$ be the operad for commutative algebras and let $M\in C(k)$.
Put $\cO=\Com\langle M,0\rangle=\Com\coprod\mathbb{F}(M,0)$. 
The algebras over $\cO$ are the commutative
algebras $A$ endowed with a map of complexes $M\to A$. Therefore, these
are the same as commutative $S(M)$-algebras and $\cO(0)=S(M)$.

We see that outside of the characteristic zero case joining a contractible
complex $M$ to an operad does not necessarily give a quasi-isomorphic
operad.

\section{}

Fortunately, the counterexample presented above is in a sense the worst 
possible case. More precisely, the following fact takes place
instead of~\cite{haha}, 6.1.1.

\begin{prp}
\label{corrected}
Let $\cO\in\Op(C(k))$. Let $M=\Cone(\id_k)[s]$
and let $n$ be a non-negative integer. Then the natural inclusion
of operads
$$ \cO\to \cO\coprod F(M,n)$$
is a componentwise quasi-isomorphism, provided either

(i) $\cO(0)=0$ and $n>0$ 

or

(ii) $k\supseteq\mathbb{Q}$.
\end{prp}

This implies the following results on model structure for $\Op(C(k))$.

\begin{thm}
\label{charzero} 
Let $k\supseteq\mathbb{Q}$. Then the category of operads in $C(k)$
admits a model structure with componentwise quasi-isomorphisms as weak 
equivalences and componentwise epimorphisms as fibrations.
\end{thm}

\begin{thm}
\label{charany} 
Let $k$ be arbitrary commutative ring. The category of operads in $C(k)$
{\em satisfying $\cO(0)=0$}
admits a model structure with componentwise quasi-isomorphisms as weak 
equivalences and componentwise epimorphisms as fibrations.
\end{thm}

\begin{rem}
A very powerful method of proving the existence of model structure on the 
category of operad is suggested by C.~Berger and I.~Moerdijk in~\cite{bm}. 

Their method allows one to prove the existence of such structure for the 
operads with $\cO(0)=0$ in any monoidal tensor category $\cE$
satisfying the following properties.
\begin{itemize}
\item $\cE$ is cofibrantly generated and its unit is cofibrant.
\item $\cE$ has a symmetric monoidal fibrant replacement functor.
\item $\cE$ admits a commutative Hopf interval.
\end{itemize}
All conditions are fulfilled if $\cE$ is the category of complexes,
so this result is more general than our Theorem~\ref{charany}.

Similarly, in the case $\cE$ admits a commutative and cocommutative Hopf 
interval, C.~Berger and I.~Moerdijk are able to prove the existence of 
model category structure as above on the category of all operads in $\cE$.

Unfortunately, the case $\cO(0)=0$ is not mentioned in~\cite{bm}. 

Instead of that, C.~Berger and I.~Moerdijk work with what they call
{\em reduced operads.} A reduced operad should mean
an operad $\cO$ endowed with an operad map $\alpha:\cO\rTo\Com$ 
such that $\alpha(0)$ is an isomorphism.

\end{rem}

\section{Proof of~Proposition~\ref{corrected}}
We will describe the operad $\cO\coprod F(M,n)$ using trees.

Here a tree is a finite oriented graph $T=(V,A,s,t)$ with the set of vertices
$V$, the set of arrows $A$, the source and the target maps $s,t:A\to V$,
such that the conditions listed below are satisfied.
\begin{itemize}
\item The graph $T$ is connected.
\item there exists a unique vertex $v_0$ (root) with no outgoing arrows.
\item Each $v\ne v_0$ admits a unique outgoing arrow.
\end{itemize}

We denote $|v|=\# s^{-1}(v)$. The set of initial vertices is
$$ \Init(T)=\{v\in V|\ |v|=0\}.$$

An {\em  $(r,n)$-marking} of the tree $T$ consists of an injective map 
$\arg:\{1,\ldots,r\}\to\Init(T)$ with the image denoted $\Arg(T)$
(this will correspond to an ordering of the arguments in an $r$-ary 
operation), and a subset $S$ of the set of vertices having valence $|v|=n$. 
The vertices from $S$ describe the slots for the elements of $M$.
Note that the vertices in $\Init(T)\setminus\Arg(T)$ are the slots
for $0$-ary operations.

A marked tree $T$ is called {\em reduced} if for each arrow $a\in A$
either $s(a)\in S$ or $s(a)\in \Arg(T)$ or $t(a)\in S$.

Let $(T,\arg,S)$ be a marked reduced tree and let $\Aut(T,\arg,S)$
(or simply $\Aut(T)$) denote the automorphism group of $(T,\arg,S)$.

We define
\begin{equation}
\label{Tcomponent}
\cO\langle M,n\rangle_T=\bigotimes_{v\in V\setminus S\setminus\Arg(T)}
\cO(|v|)\otimes M^{\otimes |S|}/\Aut(T).
\end{equation}

We claim that $\cO\coprod F(M,n)(r) $ is the direct sum of
$\cO\langle M,n\rangle_T$ over all (isomorphism classes of) $(r,n)$-marked
reduced trees.

It remains to note that if $\cO(0)=0$ and $n>0$, reduced marked trees have 
no automorphisms, and if $k\supseteq\mathbb{Q}$, all components corresponding
to non-empty $S$, are contractible.

\section{}

We know two instances where Theorem 6.1.1 from ~\cite{haha} was used
in an essential way.
In the first, \cite{tam}, 4.2.1, a version of Homology perturbation lemma
is given. Since the condition $\Char(k)=0$ is assumed throughout 
\cite{tam}, the result remains true.

Another usage of Theorem 6.1.1 from \cite{haha} is in the application
of the main theorem of \cite{voa} to realization of homotopy types.

The main idea in ~\cite{voa}, Sect.~5 was to choose a cofibrant resolution
of the Eilenberg-Zilber $\cZ$ instead of Mandell's $E_\infty$-operad
to make sure the corresponding category of algebras admits a model structure.

As we noticed in~\ref{half}, any map of operads $\cO\to \cO'$
can be factored as
$$ \cO\rTo^{i} P\rTo^{p}\cO'$$
where $i$ is a cofibration and $p$ is an acyclic fibration (i.e. surjective
quasi-isomorphism). This is enough for the existence of cofibrant
resolutions, as well as for the proof of Lemma 5.2.3 of~\cite{voa}.

Thus, the results of ~\cite{voa} do not depend on the model structure on
the category of operads.

\section{Acknowledgements}

I am grateful to B. Fresse and C. Berger for detecting an error in
Theorem 6.1.1 of ~\cite{haha}. I am grateful to I. Moerdijk who
explained me the notion of reduced operad. I am also grateful to M.~Mandell 
who pointed out that I was over-pessimistic in the first version
of the Erratum with respect to the existence of  cofibrant resolution
of the Eilenberg-Zilber operad.

\end{document}